\documentclass{amsart}

\usepackage{amsmath}
\usepackage{a4wide}

\usepackage{amsthm}
\usepackage{amsfonts}
\usepackage{eucal}
\usepackage{xspace}
\usepackage[colorlinks=true,citecolor=red,linkcolor=blue,%
urlcolor=RubineRed,pdfpagetransition=Blinds,pdftoolbar=false,pdfmenubar=false]{hyperref}

\newcommand{\field}[1]{\ensuremath{\mathbb{#1}}}

\newcommand{\C}{\field{C}\xspace}

\newcommand{\R}{\field{R}\xspace}
\newcommand{\N}{\field{N}\xspace}

\newcommand{\ens}[1]{ \left\{#1\right\} }
\newcommand{\abs}[1]{\left|#1\right|}
\newcommand{\norm}[2]{ {\left\|#1\right\|}_{#2} }

\newcommand{\ds}{\displaystyle}

\newcommand{\calL}{\mathcal L}
\newcommand{\calN}{\mathcal N}

\newtheorem{theorem}{Theorem}[section]

\newtheorem{lemma}[theorem]{Lemma}
\newtheorem{proposition}[theorem]{Proposition}

\theoremstyle{definition}
\newtheorem{remark}{Remark}[section]
\newtheorem*{question}{Question}

\begin{document}

\title[Sharp estimates for the boundary control cost]{Sharp estimates of the one-dimensional boundary control cost for parabolic systems}
\author{Manuel {\sc Gonz\'alez-Burgos}}
\thanks{Dpto.~Ecuaciones Diferenciales y An\'alisis Num\'erico and Instituto de Matem\'aticas de la Universidad de Sevilla (IMUS), Facultad de Matem\'aticas, Universidad de Sevilla, C/ Tarfia S/N, 41012 Sevilla, Spain. Supported by grant MTM2016-76990-P, Ministry of Economy and Competitiveness
(Spain). E-mail: \texttt{manoloburgos@us.es}}
\author{Lydia {\sc Ouaili}}
\thanks{Aix-Marseille Universit\'e, CNRS, Centrale Marseille, I2M, UMR 7373, 13453 Marseille, France.}

\begin{abstract}
In this work we present new results on the cost of the boundary controllability of parabolic systems at time $T > 0$. In particular, we will study optimal estimates of the control cost at time $T$ ($T$ small enough) when the eigenvalues of the generator of the $C_0$ semigroup accumulate and do not satisfy a gap condition. The main ingredient we will use is the moment method combined with sharp estimates of the $L^2(0,T; \C)$-norm of the elements of biorthogonal families to complex exponentials.
\end{abstract}

\maketitle


\section{Introduction}\label{s1}
In the last years, the controllability of non scalar systems with less controls than equations has attracted the attention of many authors. In the particular case of coupled parabolic systems, one of the main problems is obtaining necessary and sufficient conditions that allow the system to be controlled with a reduced number of distributed or boundary controls (see~\cite{AKBDGB}, \cite{FCGBdT}, \cite{AKBGBdT}, \cite{LZ},...). 

Another important problem is the study of the dependence of the so-called control cost with respect to the final observation time $T > 0$ when $T$ is small enough and the corresponding null controllability result holds at time $T > 0$. Regarding this latter problem, we highlight the works~\cite{FR1}, \cite{FR2}, \cite{Seidman}, \cite{Guichal}, \cite{Hansen}, \cite{FCZ}, \cite{Miller}, \cite{tenenbaum}, \cite{cindea}, \cite{L}, etc., where the authors study an estimate of the control cost $\mathcal{K}(T)$ in the case of scalar parabolic problems (problems that, under general assumptions, are null controllable for any $T > 0$). 
Most of the previous works uses the moment method in order to obtain an estimate of the control cost following the strategy:
\begin{enumerate}
\item Given $\ens{\lambda_k}_{k \ge 1} \subset \C$ the sequence of eigenvalues of the corresponding generator of the semigroup,  the first objective is to prove sharp bounds on the biorthogonal family to the real or complex exponentials $\Lambda = \ens{e^{- \lambda_k t}}_{k \ge 1}$ in $L^2(0,T; \C)$. To this end, in the previous works, among other assumptions, a gap condition like
	\begin{equation}\label{f15}
\inf_{k \not= n}\left| \lambda_{k} - \lambda_n \right| \ge c_0 > 0. 
	\end{equation}
is assumed.
\item The assumptions on the sequence $\Lambda$ guarantee that the scalar parabolic system under study is null controllable at any positive time $T $. Combining the moment method and the estimates on the biorthogonal family, we deduce that the following property holds: there exist positive constants $\tau_0$, $C_0$ and $C_1$ such that the control cost of the considered problem satisfies
	$$
\exp \left( \frac{C_0} T  \right) \le \mathcal{K}(T) \le \exp \left( \frac{C_1} T  \right), \quad \forall T \in (0, \tau_0) .
	$$
\end{enumerate}

The work~\cite{L} is of special relevance because in it, the author studies the cost of the controllability of the one-dimensional heat equation with a pointwise control at point $x_0$ and, in this framework, there might exist a positive minimal time of null-controllability $T_0 \in [0, \infty]$ (which depends on $x_0$ and could take any arbitrary value in $[0, \infty]$, see~\cite{Dolecki}). In this work, in particular, the author proves that, if $T_0 > 0$, the cost of the controllability at time $T > T_0$ when $T$ is close to $T_0$, may explode in an arbitrary way.

The analysis of the control cost in the framework of the controllability of coupled parabolic systems has been addressed in~\cite{BBGBO}. Again, the authors use the moment method and establish sharp bounds on biorthogonal families to complex matrix exponentials associated to complex sequences $\Lambda = \ens{\lambda_k}_{k \ge 1}$. As in the previous works, the authors impose appropriate assumptions on the sequence $\Lambda$ which include a gap condition on the terms of the sequence. These hypotheses also assure that the system under consideration is null controllable at any time $T > 0$ and, as in the scalar case, provide an estimate for $\mathcal{K}(T)$. 

However, new phenomena associated with the vectorial nature of the controllability problem arise within this framework (hyperbolic phenomena): minimal time of null controllability and dependence of the controllability result on the position of the control domain (see~\cite{MGBminT},~\cite{AKBGBdT2}, \cite{Samb},~\cite{O},...). This minimal time may come from the control action itself (as in~\cite{Dolecki}) or from the condensation index of the sequence of eigenvalues of the generator of the semigroup (see~\cite{MGBminT}). 
Let us provide more details in the case of systems with a minimal time which comes from the condensation index of the sequence. In~\cite{MGBminT}, the authors prove a general result of null controllability for abstract parabolic problems that develop a minimal time $T_0 \in [0, \infty]$ of controllability: the system is null-controllable at any time $T>T_0$ and not null-controllable for $T<T_0$. This minimal time is related to the Bernstein's condensation index of the sequence of eigenvalues $\Lambda = \ens{\lambda_k}_{k \ge 1}$ of the generator of the semigroup (see~\cite{MGBminT} and~\cite{BBM} for further details). In particular, the sequence $\Lambda $ does not satisfy, in general, the previous gap condition~\eqref{f15} and obtaining appropriate estimates of the control cost $\mathcal{K} (T) $ for $T > T_0$ is an open problem.

In this work we have two main objectives:
\begin{enumerate}
\item First, we want to give a short overview on bounds on biorthogonal families to real or complex exponentials $\ens{e^{- \lambda_k t}}_{k \ge 1}$ in $L^2(0,T; \C)$ where the sequence $\Lambda = \ens{\lambda_k}_{k \ge 1} \subset \C$ satisfies, among other hypotheses, a gap condition such as~\eqref{f15}. Let us remark that the estimate of the norm of the elements of the biorthogonal family is crucial in order to obtain the corresponding estimate of the control cost for the system under study.
\item We also want to advance new results which provide partial answers to the following questions: Given a complex sequence $ \ens{\lambda_k}_{k \ge 1}$ satisfying appropriate assumptions (see~\eqref{hyp 1}--\eqref{hyp 5},~\eqref{hyp 7} and~\eqref{hyp_6}, for an integer $q \ge 1$, below), 
\begin{enumerate}
\item Is there a biorthogonal family $\ens{q_k}_{k \ge 0} $ to $\ens{e^{- \lambda_k t}}_{k \ge 1} $ in $L^2 (0,T ; \C)$ satisfying an appropriate estimate for $\norm{q_k}{L^2 (0,T ; \C)}$ (see Theorem~\ref{t5})?
\item Using the moment method, is it possible to obtain estimates (from above and from below) of the control cost $\mathcal{K} (T)$, for $T > T_0$? ($T_0 \in [0, \infty]$ is the minimal time associated to the system).
\end{enumerate}

Observe that the kind of sequences under consideration satisfies a gap condition weaker than~\eqref{f15}. Indeed, from~\eqref{hyp_6} we deduce the existence of a positive constant $c$ (depending on $\rho$ and $q$) such that
	\begin{equation}\label{f16}
\left| \lambda_{k} - \lambda_n \right|  \ge c > 0, \quad \forall k, n: \abs{k-n} \geq q.
	\end{equation}
In particular, the sequences we are considering in this work could satisfy 
	$$
\inf_{k \not= n}\left| \lambda_{k} - \lambda_n \right| = 0.
	$$
As said before, when this happens, the null controllability problem for the system in study can have a minimal time $T_0$. 

A similar weak gap condition as~\eqref{f16} has been also considered in~\cite{BBM}. With this assumption the authors develop a block moment method to handle spectral condensation phenomenon in parabolic control problems.

\end{enumerate}

The plan of the paper is the following: In Section~\ref{s2}, we will address the null controllability problem for the one-dimensional heat equation and we present the moment method. In Section~\ref{s3} we analyze a boundary controllability problem for a non-scalar parabolic system with complex eigenvalues. For this system, we obtain the same estimate for the control cost $\mathcal{K} (T)$ as in the scalar case. In Section~\ref{s4} we study the controllability problem for the one-dimensional phase-field system with boundary controls. This is a simple example for which the eigenvalues of the corresponding generator do not satisfy the gap condition~\eqref{f15}. Section~\ref{s5} is devoted to establishing the main result of this work: sharp estimates of the norm of biorthogonal families to complex exponentials without gap conditions. Finally, in Section~\ref{s6} we apply the main result to the phase-field system and to a non-scalar parabolic system whose generator does not satisfy~\eqref{f15}. In both cases, estimates for the control cost $\mathcal{K} (T)$ are established.

%
%
%

\section{The one-dimensional heat equation: the moment method}\label{s2}
Let us first present the control cost in the simplest case: the boundary null controllability problem for the one-dimensional heat equation
	\begin{equation}\label{f1}
	\left\{ 
	\begin{array}{ll}
y_t - y_{xx} = 0  & \hbox{ in } Q_{T} = (0, \pi) \times (0,T),  \\ 
	\noalign{\smallskip}
y (0, \cdot )= v, \quad y (\pi, \cdot )=0 & \hbox{ on }  ( 0,T ), \\ 
	\noalign{\smallskip}
y (\cdot ,0 ) =y_{0}  & \hbox{ in }  ( 0,\pi ),
	\end{array}
	\right.  
	\end{equation}
where $y_0 \in H^{-1}(0, \pi)$ is given and $v \in L^{2} ( 0,T ) $ is the control function. It is well-known that system~\eqref{f1} is well-posed and the solution $y \in L^2 (Q_T)$ (defined by transposition) depends continuously on the data $y_0$ and $v$. In fact, the solution $y$ of~\eqref{f1} is related to the solutions of the corresponding adjoint problem to~\eqref{f1}
	\begin{equation}\label{f2}
	\left\{ 
	\begin{array}{ll}
-\varphi_t - \varphi_{xx} = 0  & \hbox{ in } Q_{T},  \\ 
	\noalign{\smallskip}
\varphi (0, \cdot )= \varphi (\pi, \cdot )=0 & \hbox{ on }  ( 0,T ), \\ 
	\noalign{\smallskip}
\varphi (\cdot ,T ) =\varphi_{T}  & \hbox{ in }  ( 0,\pi ),
	\end{array}
	\right.  
	\end{equation}
as follows:
%
%
\begin{theorem}\label{t1}
Let us consider problems~\eqref{f1} and~\eqref{f2}. Then, for any $y_0 \in H^{-1} (0, \pi)$ and $\varphi_T \in H_0^{1} (0, \pi)$, one has
	\begin{equation}\label{f3}
\int_0^{T} \varphi_x (0,t) v (t) \, dt  = \langle y(\cdot,T), \varphi_T \rangle - \langle y_0, \varphi(\cdot,0) \rangle ,
	\end{equation}
where $y \in L^2(Q_T) \cap C^0( [0,T] ; H^{-1} (0, \pi) )$ and  
	$$
\varphi \in L^2(0,T; H^2(0, \pi) \cap H_0^1(0, \pi)) \cap C^0( [0,T] ; H_0^1(0, \pi))
	$$ 
are, resp., the solutions to~\eqref{f1} and~\eqref{f2} associated to $y_0$ and $v$, and $\varphi_T$, and $\langle \cdot, \cdot \rangle$ stands for the duality product between $H^{-1} (0, \pi)$ and $H_0^1 (0, \pi)$.
\end{theorem}


The null controllability result for system~\eqref{f1} was obtained in~\cite{FR1} using the moment method. This method has been successfully used to prove the boundary null controllability result for general one-dimensional scalar parabolic equations with coefficients independent of $t$. Let us briefly recall this method in the case of system~\eqref{f1}.

It is well-known that the operator $- \partial _{xx}$ on $(0, \pi)$ with homogenous Dirichlet boundary conditions admits a sequence of eigenvalues and normalized eigenfunctions given by 
	\begin{equation}  \label{autov}
\lambda_{k} = k^{2}, \quad \phi _{k}(x) = \sqrt{\frac{2}{\pi}}\sin k x, \quad k\ge 1, \quad x \in (0,\pi)
	\end{equation}
which is a Hilbert basis of $L^{2} (0, \pi)$. 

From formula~\eqref{f3}, we deduce that, given $y_0 \in H^{-1}(0, \pi)$, there exists a control $v \in L^2(0,T)$ such that the solution $y$ to~\eqref{f1} satisfies $y(\cdot , T) = 0$ in $(0, \pi)$ if and only if  there exists $v\in L^2(0,T)$ satisfying 
	\begin{equation*}
\int_0^{T} \varphi_x (0,t) v (t) \, dt  = - \langle y_0, \varphi(\cdot,0) \rangle , \quad \forall \varphi_T \in H_0^1 (0, \pi),
	\end{equation*}
($ \varphi \in L^2(0,T; H^2(0, \pi) \cap H_0^1(0, \pi)) \cap C^0( [0,T] ; H_0^1(0, \pi))$ is the solution of~\eqref{f2} associated to~$\varphi_T$). Using that $\{ \phi_k \}_{k \ge 1}$ (see~\eqref{autov}) is an orthogonal basis of $H_0^1 (0, \pi)$, the previous property can be rewritten as
	\begin{equation*}
\int_0^T v(t) e^{-\lambda_k (T - t)} \partial_x\phi_k (0)\,dt = -\langle y_{0} , e^{-\lambda_k T} \phi_k \rangle , \quad \forall k \ge 1.
	\end{equation*}
In the previous equality $\lambda_k$ and $\phi_k$ are given in~\eqref{autov}. Using the Fourier decomposition of $y_{0}$, $y_{0} = \sum_{k\ge 1} y_{0,k} \phi_k$, we can conclude that there exists a control $v \in L^2(0,T)$ such that the solution $y$ to~\eqref{f1} satisfies $y(\cdot , T) = 0$ in $(0, \pi)$ if and only if  there exists $v\in L^2(0,T)$ such that
	\begin{equation*}
k \sqrt{\frac{2}{\pi}} \int_0^T e^{- \lambda_k (T-t)} v(t)\, dt = - e^{-\lambda_k T} y_{0,k}\quad \forall k \ge 1.
	\end{equation*}
Introducing the new function $u (t) = v( T-t)$, $t \in (0,T)$, the previous equality may be written as
	\begin{equation}\label{mom1}
\int_{0}^{T}e^{- \lambda_k t} u ( t ) \, dt = e^{-\lambda_k T} m_{k}, \quad \forall k\geq 1, \quad \left( m_k = - \frac 1k \sqrt{\frac{\pi}{2}}  y_{0,k} \right).
	\end{equation}
In conclusion, we have reformulated the null controllability problem for the heat equation~\eqref{f1} as a \textit{moment problem} which involves the real exponential family $\ens{ e^{- \lambda_k t} }_{k \ge 1}$. 

Let us now solve the moment problem~\eqref{mom1} in $L^2 (0,T)$, for every time $T>0$. Using that the series 
	$$
\sum_{k \ge 1} \frac{1}{\lambda_k}
	$$
is convergent, from the M\"untz Theorem (see~\cite{Schwartz}), we can conclude that this family is not total (or not complete) and minimal\footnote{We say that a sequence $\{x_n\}$ in a Banach space $X$ is minimal if no vector $x_m$ lies in the closed span of the other vectors $x_n$, i.e., for every $m \in \N$, $x_m \not\in  \overline{\hbox{span}}\, \{ x_n \}_{n \not= m}$.} in $L^2(0,T)$. This last property also characterizes the existence of, at least, a biorthogonal family $\ens{ q_k }_{k \ge 1}$ to $\ens{  e^{- \lambda_k t} }_{k \ge 1}$ in $L^2 (0,T)$ (see Lemma~5.4 in~\cite{Heil}), i.e., a family $\ens{ q_k }_{k \ge 1} \subset L^2 (0,T)$ such that
	$$
\int_0^T q_k (t) e^{- \lambda_j t} \, dt = \delta_{kj}, \quad \forall k,j \in \N.
	$$
This property provides a formal solution to the moment problem~\eqref{mom1}:
	\begin{equation}\label{f6}
u(t) = v ( T - t ) =\sum_{k\geq 1} e^{-\lambda_k T} m_{k} q_{k} ( t ) .
	\end{equation}
The question then is whether the constructed control $v$ belongs to $L^{2} ( 0,T )$. Of course, this strongly depends on the family $\ens{ q_k }_{k \ge 1}$. We have the following result:
%
%

\begin{lemma}\label{l1}
There is a positive constant $C$ such that, for all $ T > 0$ we can find a biorthogonal family $\ens{ q_k }_{k \ge 1} \subset L^2(0,T)$ to $\ens{ e^{- \lambda_k t} }_{k \ge 1}$ ($\lambda_k = k^2$) such that
	\begin{equation}\label{f4}
\| q_k \|_{L^2 (0,T)} \le C e^{C \sqrt{ \lambda_k}} e^{C/T}, \quad \forall k \ge 1.
	\end{equation}
\end{lemma}
The previous result has been proved for positive real sequences $\ens{ \lambda_k }_{k \ge 1}$ more general than~\eqref{autov}; for instance, see~\cite{FR1}, \cite{Seidman}, \cite{Miller}, \cite{tenenbaum}, ...

Throughout this work $C>0$ denotes a generic positive constant that may change line to line but which does not depend on $T$,  $k$, or $y_0$.

As a consequence of Lemma~\ref{l1} we obtain the null controllability result for system~\eqref{f1}. Indeed, from~\eqref{f4} and using the inequality
	$$ 
2 C \sqrt{\lambda_ k} \le \frac{C^2}{T} + \lambda_k T, \quad \forall C>0 \hbox{ and } k \ge 1,
	$$
we deduce that the series in~\eqref{f6} is absolutely convergent in $L^2 (0, T)$:
	$$
	\begin{array}{l}
\displaystyle \sum_{k \ge 1} \sqrt{\frac{\pi}{2}} \frac 1k  \left| y_{0,k} \right| e^{- \lambda_k T} \| q_k \|_{L^2 (0,T)} \le C \left( \sum_{k \ge 1}  \frac{|y_{0k}|^2}{k^2} \right)^{1/2} \left( \sum_{k \ge 1} e^{- 2 \lambda_k T} \| q_k \|_{L^2 (0,T)}^2 \right)^{1/2} \\
	\noalign{\smallskip}
\displaystyle \phantom{ \sum_{k \ge 1} \sqrt{\frac{\pi}{2}} \frac 1k  \left| y_{0,k} \right| } \le C e^{C /T } \| y_ 0\|_{H^{-1} (0, \pi)} \left( \sum_{k \ge 1} e^{- 2 \lambda_k T} e^{2C \sqrt{\lambda_ k}} \right)^{1/2} \\
	\noalign{\smallskip}
\displaystyle \phantom{ \sum_{k \ge 1} \sqrt{\frac{\pi}{2}} \frac 1k  \left| y_{0,k} \right| } \le C e^{C /T} \| y_ 0\|_{H^{-1} (0, \pi)}  \left( \sum_{k \ge 1} e^{- \lambda_k T} \right)^{1/2} \le \frac C{\sqrt{T}} e^{C /T} \| y_ 0\|_{H^{-1} (0, \pi)}.
	\end{array}
	$$
Therefore, the control function $v$ given in~\eqref{f6} belongs to $ L^2 (0,T) $, satisfies  
	\begin{equation}\label{f5}
\| v \|_{L^2 (0,T) } \le C e^{C /T} \| y_ 0\|_{H^{-1} (0, \pi)},
	\end{equation}
and solves the moment problem~\eqref{mom1}. We have then proved:

%
%

%
\begin{theorem}\label{t2} 
For any $y_0\in H^{-1}(0,\pi)$ and $T>0$, there exists $v\in L^{2}(0,T)$ solution to the moment problem~\eqref{mom1}. That is, $v(t)=u(T-t)$ is a null control for equation~\eqref{f1}. In addition, there exists a positive constant $C$ (independent of $T$ and $y_0$) such that~\eqref{f5} holds.
\end{theorem}
%
%


Theorem~\ref{t2} guarantees that system~\eqref{f1} can be driven from an initial datum $y_0 \in H^{-1} (0, \pi)$ to the equilibrium $0$ at time $T>0$, for any $T>0$. In particular, the set
	$$
\mathcal{Z}_T(y_0) := \left\{ v \in L^2(0, T) :  y (\cdot, T) = 0 \hbox{ in } (0,\pi), \hbox{ with } y \hbox{ solution of}~\eqref{f1} \right\},
	$$
is nonempty for any $y_0 \in H^{-1} (0, \pi)$ and any $T>0$. We can then define the control cost for system~\eqref{f1} at time $T$ as
	\begin{equation*}
\mathcal{K}(T) = \sup_{\|y_0\|_{H^{-1} (0,\pi)}=1} \left( \inf_{v\in\mathcal{Z}_{T}(y_0)} \| v \|_{ L^2(0,T)} \right), \quad \forall T>0.
	\end{equation*}

Theorem~\ref{t2} also provides an estimate for $ \mathcal{K}(T) $ and a measure of how violent fast controls for system~\eqref{f1} are:
	$$
\mathcal{K}(T) \le C e^{C /T} , \quad \forall T > 0.
	$$

The previous estimate for the control cost in $L^2(0,T)$ is optimal for problem~\eqref{f1}. Indeed, G\"uichal proved in~\cite{Guichal}, the optimality of the previous expression with respect to the time dependence:
\begin{theorem}[G\"uichal] \label{t3}
Under the previous conditions, there exist positive constants $C_0$ and $C_1$ independent of $T > 0$ such that, for any $T \in (0, 1]$ one has
	\begin{equation}\label{f7}
e^{C_0 / T} \le \mathcal{K}(T) \le e^{C_1 /T}, \quad \forall T \in (0, 1].
	\end{equation}
\end{theorem}

\begin{remark}  
The estimate~\eqref{f4} can be obtained for more general positive real sequences $\{ \lambda_k \}_{k \ge 1}$. To be precise, in order to obtain~\eqref{f4} we have used:
\begin{enumerate}
\item $\{ \lambda_k \}_{k \ge 1} \subset (0, \infty)$ is an increasing sequence and there exist $K > 0$ and $ \alpha \in \R $ such that
	\begin{equation}\label{raiz}
\sqrt{\left| \lambda_k \right| } = K (k + \alpha)\left[ 1+  O(1/k) \right]. 
	\end{equation}
\item Gap condition: There is $\rho > 0$ such that
	\begin{equation}\label{gap}
\left| \lambda_k - \lambda_n \right| \ge \rho \left| k^2 - n^2 \right|, \quad \forall k,n \ge 1. 
	\end{equation}
\end{enumerate}
In particular, Lemma~\ref{l1}, estimate~\eqref{f4} and Theorem~\ref{t2} can be also proved for the corresponding sequences of eigenvalues of  more general second order elliptic operators: 
	$$
L y = -\left( p( \cdot ) y_x \right)_x + q( \cdot ) y, \quad D(L) = H^2 (0, \pi) \cap H_0^1(0, L),
	$$
with $p \in W^{1, \infty} (0, \pi)$ and $q \in L^\infty (0, \pi)$ satisfying $p (\cdot) \ge p_0 > 0$ in $(0, \pi )$.
\end{remark}

\section{A boundary controllability problem for a parabolic system with complex eigenvalues}\label{s3}
In order to extend inequalities~\eqref{f4} and~\eqref{f7} to complex sequences or sequences that do not satisfy conditions~\eqref{raiz} or~\eqref{gap}, let us consider the following boundary controllability problem:
	\begin{equation}\label{f8}
\left\{
	\begin{array}{ll}
y_t - y_{xx} + A_0 y = 0 & \hbox{in } Q_T = (0, \pi) \times (0,T), \\ 
	\noalign{\smallskip}
y (0, \cdot) = B v , \quad  y (\pi , \cdot) = 0 & \hbox{on } ( 0,T) , \\ 
	\noalign{\smallskip}
y( \cdot , 0) =y_{0}, & \hbox{in }( 0,\pi ) .
	\end{array}
\right.
	\end{equation}
In~\eqref{f8} $ y = (y_1, y_2)^t$, $y_0 \in H^{-1}((0,\pi); \R^2)$, $v \in L^2(0,T)$ and 
$$
	A_0 =  \left( 
	\begin{array}{cc} 
0 & 1 \\ -1 & 0  
	\end{array}\right), \ B = \left( 
	\begin{array}{c}
 0 \\ 1  
	\end{array}\right).
$$
As before, we want to prove a null controllability result for system~\eqref{f8} but, unlike system~\eqref{f1}, we only exert a scalar control on the second equation while the first equation is indirectly controled by means of the coupling matrix $A_0$. 

Again, system~\eqref{f8} is well-posed. One has: 
\begin{proposition}\label{p1}
Assume $y_0 \in H^{-1}(0,\pi;\R^2)$ and $v \in L^2(0,T)$. Then, system~\eqref{f8} admits a unique solution (by transposition) $ y \in L^2(Q_T;\R^2) \cap C^0([0,T]; H^{-1}(0,\pi; \R^2))$ which depends continuously on the data. 
\end{proposition}

We can repeat the arguments of Section~\ref{s2} and prove that the null controllability of system~\eqref{f8} at time $T>0$ is equivalent to an appropriate moment problem. In this case, the adjoint operator $L^* $ is given by
	$$
L^* \varphi = - \varphi_{xx} + A_0^* \varphi , \quad D(L^* ) = H^2 (0, \pi; \C^2) \cap H_0^1(0, \pi;  \C^2), 
	$$
and admits a complex sequence of eigenvalues and normalized eigenfunctions given by 
	\begin{equation*}
\Lambda = \{ n^{2}- i , n^2+i \}_{n \ge 1} \quad \mathcal{B} = \ens{ \frac{1}{\sqrt{\pi}}\sin n x  \left( \begin{array}{c} i \\ -1 \end{array} \right), \frac{1}{\sqrt{\pi}}\sin n x \left( \begin{array}{c} 1 \\ - i \end{array} \right) }_{n \ge 1} .
	\end{equation*}
In addition, $\mathcal{B} $ is a Riesz basis of $L^{2} (0, \pi; \C^2)$. 

Let us take 
	$$
\lambda_k = \left\{
	\begin{array}{ll}
\displaystyle \left( \frac{k+1}{2} \right)^2 - i , & \hbox{if $k$ is odd},   \\
	\noalign{\smallskip}
\displaystyle \left( \frac{k}{2} \right)^2 + i , & \hbox{if $k$ is even.}
 	\end{array}
\right. 
	$$
Then, it is not difficult to see that the sequence $\Lambda$ of eigenvalues of $L^*$ is given by $\Lambda = \ens{ \lambda_{k} }_{k \ge 1} \subset \C $ and satisfies~\eqref{raiz} for $K = 1/2$ and $\alpha = 0$. Nevertheless, the sequence $\Lambda $ does not satisfy the gap condition~\eqref{gap}. 

Now the question is: Would it be possible to obtain Lemma~\ref{l1} for the previous complex sequence $\Lambda$ in $L^2(0,T; \C)$? In this regard, one has the following generalization of Lemma~\ref{l1}:
\begin{lemma}\label{l2}
Let $\Lambda = \ens{\lambda_k}_{k \geq 1} \subset \C$ be a sequence of complex numbers with the following properties:
\begin{enumerate}
\item\label{hyp 1}
$\lambda_k \neq \lambda_n$ for all $k,n\in \N$ with $k \neq n$.
\item\label{hyp 2}
$\Re(\lambda_k) >0$ for every $k \geq 1$.

\item\label{hyp 4}
For some $\beta>0$,
	$$
\abs{\Im (\lambda_k ) } \le \beta \sqrt{\Re (\lambda_k )}, \quad \forall k \ge 1.
	$$
\item\label{hyp 5}
$\ens{\lambda_k}_{k \geq 1}$ is non-decreasing in modulus: 
	$$
\abs{\lambda _k} \leq \abs{\lambda_{k+1}}, \quad \forall k\geq 1.
	$$
\item\label{hyp 6}
$\ens{\lambda_k}_{k \geq 1}$ satisfies the following gap condition: for some $\rho,q>0$,
	$$
	\left\{
	\begin{array}{l}
\ds \abs{\lambda_k - \lambda_n} \ge \rho \abs{k^2 - n^2}, \quad \forall k, n: \abs{k-n} \geq q. \\
	\noalign{\smallskip}
\ds \inf_{k \neq n:\abs{k-n}<q} \abs{\lambda_k-\lambda_n}>0.
	\end{array}\right.
	$$
\item\label{hyp 7}
For some $p,\alpha>0$,
	\begin{equation}\label{counting-property}
\abs{p \sqrt{r} - \calN (r)} \le \alpha, \quad \forall r > 0,
	\end{equation}
where $\calN$ is the counting function associated with the sequence $\ens{\lambda_k}_{k \geq 1}$, that is the function defined by
	\begin{equation}\label{counting}
\calN (r) = \# \ens{k : \abs{\lambda_k} \le r} , \quad \forall r >0.
	\end{equation}
\end{enumerate}
Then, there exists $\tau_0>0$ such that, for every $0<T<\tau_0$, we can find a family of $\C$-valued functions
	$$
\ens{ q_{k} }_{k\geq 1 }\subset L^2(0 , T; \C)
	$$ 
biorthogonal to $\ens{e^{- \lambda_k t }}_{k\geq 1}$ in $L^2(0 , T; \C)$, with in addition
	\begin{equation} \label{biacot}
\norm{q_{k}}{L^2(0, T; \C)} \le C e^{ C \sqrt {\Re( \Lambda_k )}  + \frac{C}{T} }, \quad \forall k \ge 1.
	\end{equation}
\end{lemma}

The proof of this result can be found in~\cite{BBGBO}. 

We can now apply Lemma~\ref{l2} to problem~\eqref{f8} in order to obtain the null controllability result for this system at time $T > 0$ and an estimate of the control cost at this time. Indeed, the sequence $\Lambda$ of eigenvalues of $L^* = - \partial_{xx} + A_0 $ satisfies assumptions~\eqref{hyp 1}--\eqref{hyp 7} in~Lemma~\ref{l2}. On the other hand, as said before, the null controllability at time $T >0$ of system~\eqref{f8} is equivalent to a moment problem for the complex exponentials $\ens{e^{- \lambda_k t} }_{k \ge 1}$ in $L^2(0,T ; \C)$. If we now combine inequality~\eqref{biacot} and the reasoning of Section~\ref{s2}, we deduce:

%
%

%
\begin{theorem}\label{t4} 
System~\eqref{f8} is null controllable at time $ T > 0 $, for any $T>0$. In addiction, there exists a constant $C >0$ such that the control cost for system~\eqref{f8} satisfies 
	\begin{equation} \label{f9}
\mathcal{K} (T) \le C e^{C/T}, \quad \forall T > 0. 
	\end{equation}
\end{theorem}
\begin{remark}\label{r1}
The null controllability result for system~\eqref{f8} at time $T>0$ was established in~\cite{FCGBdT}. To be precise, in this reference the authors provided a necessary and sufficient condition for the boundary null controllability of general $2 \times 2$ parabolic systems with constant coefficients: 
	\begin{equation}\label{f10}
\left\{
	\begin{array}{ll}
y_t - y_{xx} + A y = 0 & \hbox{in } Q_T = (0, \pi) \times (0,T), \\ 
	\noalign{\smallskip}
y (0, \cdot) = B v , \quad  y (\pi , \cdot) = 0 & \hbox{on } ( 0,T) , \\ 
	\noalign{\smallskip}
y( \cdot , 0) =y_{0}, & \hbox{in }( 0,\pi ) ,
	\end{array}
\right.
	\end{equation}
with $ A \in \calL(\R^2)$, $B \in \R^2$, $y_0 \in H^{-1}(0,\pi; \R^2)$ and $v \in L^2(0,T)$. This necessary and sufficient condition for the boundary null controllability of system~\eqref{f10} was generalized in~\cite{AKBGBdT} for general $n \times n$ systems with the same structure as~\eqref{f10} with $A \in \calL(\R^n) $ and $B \in \calL (\R^m ; \R^n)$ ($n,m \ge 1$). Nevertheless, in these two references the authors did not prove the estimate of the control cost~\eqref{f9}. This estimate was proved in~\cite{BBGBO}. In fact, the authors proved a generalization of Lemma~\ref{l2} for complex exponentials $\ens{t^je^{- \lambda_k t}}_{k\ge 1, 0 \le j \le \eta}$, with $\eta \ge 0$ a given integer, and, as a consequence, they showed~\eqref{f9} when system~\eqref{f10} satisfies the sufficient null controllability condition. 
\end{remark}

\begin{remark}\label{r2}
Assumption~\eqref{hyp 6} (gap condition) on the sequence $\Lambda$ is crucial in the proof of the existence of a biorthogonal family $\ens{q_k}_{k \ge 1}$ to $\ens{e^{- \lambda_k t}}_{k \ge 1}$ in $L^2(0,T; \C)$ that satisfies~\eqref{biacot}. However, we will see in the next section that there are some boundary controllability problems for physical systems which do not satisfy the previous gap condition.  Thus, it seems pertinent to analyze the control cost for systems that do not satisfy~\eqref{hyp 6}.\end{remark}

%
%
%

\section{The boundary controllability of a phase-field system}\label{s4}
In this section we will study the boundary null controllability of a phase field system of Caginalp type (see~\cite{cag}) which is a model describing the transition between the solid and liquid phases in solidification/melting processes of a material occupying the  interval $(0, \pi)$. For that purpose, we consider the nonlinear system
	\begin{equation}\label{PFSy}
	\left\{
	\begin{array}{ll}
\displaystyle {\theta}_t - \xi{\theta}_{xx} + \dfrac{1}{2}\rho\xi{\phi}_{xx} + \dfrac{\rho}{\tau}{\theta} = f({\phi})	& \mbox{in } Q_T, 
	\\
	\noalign{\smallskip}
\displaystyle {\phi}_t - \xi{\phi}_{xx} - \dfrac{2}{\tau}{\theta} = -\frac 2 \rho f({\phi}) & \mbox{in }  Q_T, 
	\\
	\noalign{\smallskip}
\displaystyle {\theta}(0,\cdot) = v,\ {\phi}(0,\cdot) = c,\ {\theta}(\pi,\cdot)=0 , \ {\phi}(\pi,\cdot) = c & \mbox{on }  (0,T),	 
	\\
	\noalign{\smallskip}
\displaystyle {\theta}(\cdot,0) = {\theta}_0, \  {\phi}(\cdot,0)={\phi}_0 & \mbox{in }  (0,\pi),	\end{array}
	\right. 
	\end{equation}
where: ${{\theta} = {\theta}(x,t)}$ is the temperature of the material; ${{\phi} = {\phi} (x,t)}$ is the phase-field function used to identify the solidification level of the material; ${c \in \{ - 1,0,1\}}$; ${f}$ is the nonlinear term which comes from the derivative of the classical regular double-well potential $W$:
	$$ 
\displaystyle f({\phi}) = -\frac{\rho}{4\tau}\left( {\phi}-{\phi}^3 \right) . 
	$$
On the other hand, 
$\rho > 0$, $ \tau > 0$ and $ \xi > 0$ are, resp., the latent heat, a relaxation time and the thermal diffusivity. Finally, ${v \in L^2(0,T)}$ is the control function, and ${ \theta_0,  \phi_0}$ are the initial data.

The null controllability property of the nonlinear system~\eqref{PFSy} has been analyzed in~\cite{GBSN} and depends on the coefficients $\rho$, $\tau$ and $ \xi $. This property can be obtained from the corresponding one of the linear version of~\eqref{PFSy} (see~\cite{GBSN} for more details):
	\begin{equation}\label{ff}
	\left\{ 
	\begin{array}{ll}
y_{t}-D y_{xx} + A y= 0  & \hbox{in }Q_{T}:= (0, \pi) \times (0,T), \\ 
	\noalign{\smallskip} 
y(0, \cdot ) = B v ,\quad y(\pi, \cdot) =0  & \hbox{on } (0,T), \\ 
	\noalign{\smallskip} 
y(\cdot ,0)=y_{0} & \hbox{in }(0, \pi) ,
	\end{array}
	\right. 
	\end{equation}
where $y = (\theta , \phi)$,
	\begin{equation}\label{ADB}
D = \left(
	\begin{array}{cc}
\xi & -\dfrac{1}{2} \rho \xi \\
	\noalign{\smallskip}
0 &  \xi
	\end{array}
	\right),
\quad
A = \left(
	\begin{array}{cc}
\dfrac{\rho}{\tau} & -\dfrac{\rho}{2\tau} \\
	\noalign{\smallskip}
-\dfrac{2}{\tau} & \dfrac{1}{\tau}
	\end{array}
	\right),
\quad
B =
	\left(
	\begin{array}{cc}
1 \\ 0
	\end{array}
	\right). 
	\end{equation} 

Again, our objective is to drive the solution $y = (\theta , \phi)$ of~\eqref{ff} to the equilibrium $(0,0)$ only acting on the first equation (the temperature) through the Dirichlet boundary condition at point $x =0$. The phase $\phi$ is indirectly controlled using the coupling matrices $D$ and $A$.

Let us first establish the well-posedness of  system~\eqref{ff}. One has:
\begin{proposition}\label{p2}
Assume $y_0 = (\theta_0, \phi_0) \in H^{-1}(0,\pi;\R^2)$ and $v \in L^2(0,T)$. Then, system~\eqref{ff} admits a unique solution (by transposition) $ y = (\theta, \phi) \in L^2(Q_T;\R^2) \cap C^0([0,T]; H^{-1}(0,\pi; \R^2) )$ which depends continuously on the data. 
\end{proposition}

Now, we can already analyze the controllability properties of system~\eqref{ff}. The approximate controllability of this system  is given by the next result:
\begin{theorem}[Approximate controllability] \label{CAproximada_}
Fix $T>0$. Then, system~\eqref{ff} is approximately controllable in $H^{-1} (0, \pi; \R^2)$ at time $T > 0$ if and only if
	\begin{equation}\label{H2}
{\xi^2\tau^2(\ell^2 - k^2)^2 - 2\xi\rho\tau(\ell^2 + k^2)- 2 \rho-1 \neq 0}, \quad \forall k,\ell\geq1, \quad \ell >k.
	\end{equation}
\end{theorem}

The proof of this result can be found in~\cite{GBSN}.

It is well-known that the approximate controllability at time $T$ for system~\eqref{ff} is necessary in order to prove the null controllability result for this system at same time $T$. So, assuming condition~\eqref{H2}, the next task will be to study the null controllability of the linear problem~\eqref{ff} at time $T > 0$ and the corresponding control cost. Again, these properties can be obtained from the properties of the sequence of eigenvalues of the operators $L$ and $L^*$ with $L := -D  \partial_{xx} + A $ and $D(A) = H^2 (0,\pi ; \R^2) \cap H_0^1 (0,\pi ; \R^2) $. 

	\begin{proposition}\label{Lspec}
Let us consider the operators $L := -D  \partial_{xx} + A$ and $L^*$ (the matrices $D$ and $A$ are given in~\eqref{ADB}). Then, the spectra of $L$ and $L^*$ are given by $\sigma(L) = \sigma (L^*)= \ens{\lambda_k^{(1)}, \lambda_k^{(2)} }_{k\geq1}$ with
	\begin{equation}\label{lamk12}
\lambda_k^{(1)} = \xi k^2 + \dfrac{\rho + 1}{2\tau} - r_k ,\quad \lambda_k^{(2)} = \xi k^2 + \dfrac{\rho + 1}{2\tau} + r_k, \quad \forall k \ge 1,			
	\end{equation}
where $r_k$ is given by
	\begin{equation*}
r_k:= \sqrt{ \dfrac{\xi\rho}{\tau} k^2 + \left(\dfrac{\rho+1}{2\tau}\right)^2 } , \quad \forall k \ge 1.
	\end{equation*}
In addition, $\lambda_k^{(1)} \not= \lambda_n^{(2)}$ for any integers $k,n \ge 1$ if and only if~\eqref{H2} holds.
\end{proposition}

Again, the proof of this result can be found in~\cite{GBSN}.

Thus, the sequence of eigenvalues of the vectorial elliptic operator $L$ is real. If condition~\eqref{H2} is fullfilled, it can be rearranged in such a way that $\Lambda = \ens{\lambda_k}_{k \ge 1} = \ens{\lambda_k^{(1)}, \lambda_k^{(2)}}_{k \ge 1} $ is an increasing sequence that satisfies~\eqref{hyp 1}--\eqref{hyp 5} and~\eqref{hyp 7} for appropriate $p, \alpha >0 $. However, in general
	$$
\ds \inf_{k \ge 1} (\lambda_{k + 1} - \lambda_k) = 0,
	$$
and condition~\eqref{hyp 6} does not hold (see~\cite{GBSN} for the details). Therefore, we have a new system where the associated sequence of eigenvalues does not satisfy the gap condition~\eqref{hyp 6} and Lemma~\ref{l2} cannot be applied. Thus, the following question arises:

\begin{question}
When the gap condition~\eqref{hyp 6} fails, is system~\eqref{ff} null controllable at time $T$ for any $T>0$?; in this case, can we obtain an estimate of the control cost? We will answer this two question in the next sections.
\end{question}

\begin{remark}[Minimal time]\label{r3}
Some choices of the parameters $\rho$, $\tau$ and $\xi$ make the eigenvalues of the operators $L := -D  \partial_{xx} + A$ and $L^*$ (the matrices $D$ and $A$ are given in~\eqref{ADB}) concentrate and this could affect the null controllability result. To be precise, in~\cite{MGBminT}, the authors proved that when the eigenvalues of the generator of a $C_0$ semigroup accumulate, the corresponding null controllability problem has a minimal time $T_0 \in [0, \infty]$ of null controllability. This hyperbolic behaviour of the parabolic system is related to the condensation index of the sequence of eigenvalues of the generator. In the case of system~\eqref{ff}, this minimal time is $T_0 = 0$ (see~\cite{GBO}).
\end{remark}

\begin{remark}\label{r4}
The previous example shows that it is necessary to generalize Lemma~\ref{l2} for complex sequences $\Lambda = \ens{\lambda_k}_{k \geq 1} \subset \C$ that satisfy all the assumptions in this lemma except the gap condition:
	$$
\inf_{k \neq n:\abs{k-n}<q} \abs{\lambda_k-\lambda_n} > 0.
	$$
For this kind of sequences, it would also be interesting to analyze the dependence of the control cost $ \mathcal{K}(T) $ with respect to time $ T> 0 $, when $ T $ is small. 
\end{remark}

%
%
%

\section{The main result: Bounds on biorthogonal families to complex exponentials without gap condition}\label{s5}
In this section we will deal with the main result of this work. To this end, let us consider a complex sequence $\Lambda = \ens{\lambda_k}_{k \ge 1} \subset \C$ satisfying~\eqref{hyp 1}--\eqref{hyp 5}, ~\eqref{hyp 7} and
	\begin{equation}\label{hyp_6}
\ds \abs{\lambda_k - \lambda_n} \ge \rho \abs{k^2 - n^2}, \quad \forall k, n: \abs{k-n} \geq q, 
	\end{equation}
with $q \ge 1$, a given integer. Our objective is to prove a result like Lemma~\ref{l2} for the sequence $\Lambda$ but without imposing the gap condition~\eqref{hyp 6}.

\begin{remark}
Observe that we do not impose assumption~\eqref{hyp 6} on the sequence $\Lambda$. Therefore, the infimum
	$$
\inf_{k \not= n} \left| \lambda_k - \lambda_n \right|
	$$
could be zero. Nevertheless, it is possible to prove that the sequence $\Lambda$ satisfies
	$$
\sum_{k \ge 1} \frac{1}{| \lambda_k|} < \infty, 
	$$
and the exponential sequence $\ens{e^{- \lambda_k t}}_{k \ge 1}$ admits a biorthogonal family $\ens{q_k}_{k \ge 1}$ in $L^2 (0, T; \C)$, for any $T >0$. This fact provides a formal solution $u$ to the moment problem 
	\begin{equation*}
\int_{0}^{T}e^{-\lambda_k t} u ( t ) \, dt = e^{-\lambda_k T} m_{k}, \quad \forall k\geq 1, 
	\end{equation*}
with $\ens{m_k}_{k \ge 1} \subset \ell^2$, given by
	$$
u(t) = \sum_{k\geq 1} e^{-\lambda_k T} m_{k} q_{k} ( t ) .
	$$
In fact, the previous expression, in general, defines a function in $L^2 (0,T ; \C)$ provided $T$ is large enough. To be precise, one can prove that there exists a minimal time $T_0 \in [0, \infty ]$ such that the corresponding problem is null controllable at time $T > 0$ if $ T > T_0$.  On the other hand, if $T < T_0$, the system is not null controllable at time $T$. For more details, see~\cite{FCGBdT},~\cite{AKBGBdT},~\cite{MGBminT}, ~\cite{BBM} and~\cite{GBO}.

Even in the case in which the minimal time is $T_0 = 0$, we cannot apply Lemma~\ref{l2} to the sequence $\Lambda$ and deduce inequality~\eqref{f9} for the corresponding control cost.
\end{remark}

One has:
\begin{theorem}\label{t5}
Assume that the complex sequence $\Lambda = \ens{\lambda_k}_{k \ge 1} \subset \C$ satisfies~\eqref{hyp 1}--\eqref{hyp 5},~\eqref{hyp 7} and~\eqref{hyp_6} for an integer $q \ge 1$. Then, there exists a constant $C >0$, only depending on $\beta$, $q$, $\rho$, $p$ and $\alpha$, such that for any $T  > 0$ there exists a biorthogonal family $\ens{q_k}_{k \ge 1} $ to $\ens{e^{- \lambda_ k t}}_{k \ge 1}$ in $L^2(0,T; \C)$ such that
	$$
\norm{q_k}{L^2(0,T; \C)} \le C e^{C  \sqrt{\Re(\lambda _{k})} } e^{C /T} \prod_{1 \le | k - n | < q} | \lambda_k - \lambda_n |^{-1} , \quad \forall k\geq 1. 
	$$
\end{theorem}
The proof is very technical and follows some ideas from~\cite{BBGBO}. It can be found in~\cite{GBO}.

\begin{remark}
Theorem~\eqref{t5} generalizes Lemma~\ref{l1}, Lemma~\ref{l2} and some previous works. For instance:
\begin{enumerate}
\item \cite{FR1}, \cite{Seidman}, \cite{Miller} and \cite{tenenbaum} where the authors considered real positive sequences satisfying~\eqref{hyp_6} for $q = 1$, i.e., for any positive integers $k,n \ge 1$, with $k \not= n$;
\item \cite{Hansen} and~\cite{Seidman00} where the authors considered complex sequences that, again, satisfy a gap condition;
\item \cite{cindea}: in this case, the authors proved a similar result to that established in Theorem~\ref{t5} but for positive real sequences satisfying the conditions
	\begin{equation*}
	\left\{
	\begin{array}{l}
\displaystyle \left| \lambda_k - \alpha k^2 \right| \le c k, \quad \forall k \ge 1, \\
  	\noalign{\smallskip}
\displaystyle \inf_{n : n \not= k} \left| \sqrt{\lambda_k} - \sqrt{\lambda_n} \right| \ge \frac rk, \quad \forall k \ge 1,  
	\end{array}
	\right.
	\end{equation*}
with $\alpha $, $c$ and $r$ positive constants independent of $k$. In particular, the previous assumptions imply that the sequence $\ens{\lambda_k}_{k \ge 1}$ does not condense and satisfies the gap condition~\eqref{f15}. Indeed, for $k \ge 2c / \alpha$ and any $n \ge 1$, one has
	$$
\lambda_k \ge \alpha k^2 - c k \ge \frac 12 \alpha k^2,
	$$
and
	$$
\left| \lambda_k - \lambda_n \right| =  \left| \sqrt{\lambda_k} - \sqrt{\lambda_n} \right|  \left( \sqrt{\lambda_k} + \sqrt{\lambda_n} \right) \ge \frac rk \sqrt{\lambda_k} \ge r \sqrt{\frac{\alpha}{2}} > 0.
	$$
On the other hand, if $k < 2c / \alpha $, inequality~\eqref{f15} can be easily obtained.
\end{enumerate}
\end{remark}

%
%
%

\section{Control cost for some parabolic systems}\label{s6}
In this section we will analyze the control cost for some parabolic systems when we exert a boundary control on the system. The estimates of the control cost are a consequence of Theorem~\ref{t5}.


\subsection{The linear phase-field system}\label{s6.1}
In this section we again consider system~\eqref{ff} with $D$, $A$ and $B$ given in~\eqref{ADB}. Now, our objective is to give a null controllability result at time $T > 0$ for this system when~\eqref{H2} holds and obtain a bound for the corresponding control cost. One has:
\begin{theorem}\label{t6}
Under the previous conditions, assume that condition~\eqref{H2} holds. Then, system~\eqref{ff} is exactly controllable to zero at any time $T>0$. Moreover, there exists a positive constant $C$ such that 
	$$
\mathcal{K} (T) \le C e^{C/T}, \quad \forall T>0,
	$$
where $\mathcal{K} (T) $ is the control cost for system~\eqref{ff}:
	\begin{equation*}
\mathcal{K}(T) = \sup_{\|y_0\|_{H^{-1} (0,\pi; \R^2 )}=1} \left( \inf_{v\in\mathcal{Z}_{T}(y_0)} \| v \|_{ L^2(0,T)} \right), \quad \forall T>0.
	\end{equation*}
and
	$$
\mathcal{Z}_T(y_0) := \left\{ v \in L^2(0, T) :  y (\cdot , T) = 0 \hbox{ in } (0,\pi), \hbox{ with } y \hbox{ solution of }~\eqref{ff} \right\}.
	$$
\end{theorem}

\begin{proof}[Sketch of the proof]
As saw in Section~\ref{s2}, given $T > 0 $ and $y_0 \in H^{-1}(0,\pi; \R^2 )$, there exists a control $v \in L^2(0,T)$ such that the solution of~\eqref{ff} satisfies $y (\cdot, T) = 0$ in $(0, \pi)$ if and only if $v \in L^2(0, T) $ solves the moment problem
	\begin{equation*}
\int_{0}^{T}e^{-\lambda_k t} u ( t ) \, dt = e^{-\lambda_k T} m_{k}, \quad \forall k\geq 1,  \quad \left( m_k = m_k(y_0) \hbox{ is s.t.} \ens{m_k}_{k \ge 1} \in \ell^2 \right),
	\end{equation*}
with $\Lambda = \ens{\lambda_k}_{k \ge 1} = \ens{\lambda_k^{(1)}, \lambda_k^{(2)}}_{k \ge 1} $ and $\lambda_k^{(i)}$ given in~\eqref{lamk12}. On the other hand, the real positive sequence $\Lambda$ satisfies conditions~\eqref{hyp 1}--\eqref{hyp 5},~\eqref{hyp 7},~\eqref{hyp_6} (for $q = 2$) and
	$$
\abs{\lambda_k - \lambda_{k-1}} \abs{\lambda_k - \lambda_{k+1}} \ge \frac Ck \quad (C = C (\xi, \rho, \tau) >0 ),
	$$
(see~\cite{GBSN} for the details). We can therefore apply Theorem~\ref{t5} to the sequence $\Lambda$ and conclude the existence of a positive constant $C$ and a biorthogonal family $\ens{q_k}_{k \ge 1}$ to $\ens{e^{- \lambda_ k t}}_{k \ge 1}$ in $L^2(0,T)$ such that
	$$
\| q_k \|_{L^2(0,T)} \le C e^{C  \sqrt{\Re(\lambda _{k})} } e^{C /T} \prod_{1 \le | k - n | < 2} | \lambda_k - \lambda_n |^{-1}  \le C e^{C  \sqrt{\Re(\lambda _{k})} } e^{C /T}, \quad \forall k\geq 1. 
	$$

Combining the previous inequality and the reasoning in Section~\ref{s2} we deduce the existence of a control $v \in L^2 (0,T)$ which solves the previous moment problem and
	$$
\norm{v}{L^2(0,T)} \le C e^{C/T} \norm{y_0}{H^{-1} (0, \pi; \R^2)} ,
	$$
with $C >0$ a positive constant only depending on $\xi$, $\rho$ and $\tau$. This provides the estimate of $\mathcal{K} (T)$ for system~\eqref{ff} and the proof of Theorem~\ref{t6}
\end{proof}

\begin{remark}
In~\cite{GBSN} the authors prove Theorem~\ref{t6} imposing~\eqref{H2} and an additional condition to the parameters~$\xi$, $\rho$ and $\tau$ that, in particular, implies
	$$
\inf_{k \not= n} \left| \lambda_k - \lambda_n \right| \ge c_0 >0
	$$
for a positive constant $c_0 > 0$. Therefore, Theorem~\ref{t6}generalizes the null controllability result for the linear system~\eqref{ff} proved in the previous reference. 
\end{remark}


\subsection{A $ 2 \times 2$ linear coupled parabolic system}\label{s6.2}
In this section we are going to give an interesting example in which the fast controls are more violent than those of the heat equation. We will see that this violent behaviour comes from the condensation of the eigenvalues of the corresponding vectorial elliptic operator.

Let us consider the problem
	\begin{equation}\label{f11}
	\left\{
	\begin{array}{ll}
y_t - y_{xx} + a(x) A_1 y = 0 & \hbox{in } (0, \pi) \times (0,T), \\
	\noalign{\smallskip}
y(0, \cdot) = B v , \quad y(\pi, \cdot) = 0  & \mbox{on } (0,T), \\
	\noalign{\smallskip}
y(\cdot ,0)=y_0 &\mbox{in } (0, \pi), 
	\end{array}
	\right.
	\end{equation}
where $a \in L^2 (0, \pi)$ is a function such that 
	\begin{equation}\label{f12}
\displaystyle \int_0^\pi a (x) \, dx =0,
	\end{equation}
$y = (y_1, y_2)^t$, $y_0 \in H^{-1}(0, \pi; \R^2)$, $v \in L^2 (0, T)$ and
	$$
A_1 =\left( 
\begin{array}{cc}
0 & 0 \\ 
0 & 1
\end{array}%
\right)  , \quad B=\left( 
\begin{array}{c}
1 \\ 
1
	\end{array}%
\right).
	$$

As before, our objective is to drive the system to zero at time $T >0$ acting on the Dirichlet boundary condition at $x = 0$ with a scalar control $v \in L^2 (0,T)$. 

System~\eqref{f11} is well-posed and for any $y_0 \in H^{-1} (0, \pi ; \R^2)$ and $v \in L^2(0,T)$, the system has a unique solution (by transposition) $ y \in L^2(Q_T;\R^2) \cap C^0([0,T]; H^{-1}(0,\pi; \R^2) )$ which depends continuously on the data. 

The null controllability of system~\eqref{f11} has been studied in~\cite{O}. In this work, the author proves that, in general, system~\eqref{f11} has a minimal time of null controllability $T_0 \in [0, \infty]$ when $a$ satisfies~\eqref{f12}. Nevertheless, 
obtaining an estimate of the control cost for system~\eqref{f11} in $L^2(0,T )$ has not been addressed before and is an open problem when~\eqref{f12} holds.

As a consequence of Theorem~\ref{t5}, let us provide an estimate of this control cost in a simple case. Before, let us summarize the controllability result for system~\eqref{f11} proved in~\cite{O}. 

Let us introduce the vectorial operator
	$$
L := - \partial_{xx} + a (\cdot) A_1 \quad  \hbox{with } D(L) = H^2(0,\pi; \R^2) \cap H_0^1(0,\pi; \R^2).
	$$
Under assumption~\eqref{f12}, one has,
	$$
\sigma(L) = \sigma (L^*) = \Lambda = \ens{k^2, k^2 + \beta_k }_{k \ge 1}, \hbox{with } \ens{\beta_k}_{k \ge 1} \in \ell^2, 
	$$
(see~\cite{O}). Since $\lim \beta_k = 0$, the sequence $\Lambda$ condensates at infinity and the gap condition~\eqref{hyp 6} does not hold. This property has strong implications in the controllability of the system. One has:

\begin{theorem}\label{t7}
Let us consider $a \in L^2(0, \pi)$, a function satisfying~\eqref{f12}. Given $T >0$, one has:
\begin{enumerate}
\item System~\eqref{f11} is approximately controllable at time $T > 0$ if and only if $ \Lambda_1 \cap \Lambda_2 = \emptyset$, where
	$$
\Lambda_1 = \sigma \left( - \partial_{xx} \right) = \ens{k^2}_{k \ge 1} \quad \hbox{and} \quad  \Lambda_2 = \sigma \left(- \partial_{xx} + a(\cdot) \right) = \ens{k^2 + \beta_k }_{k \ge 1} .
	$$
\item Assume that $ \Lambda_1 \cap \Lambda_2 = \emptyset$ and define
	$$
T_0 = \limsup \frac{- \log \left( \abs{\beta_k} \right)}{k^2} \in [0, \infty] .
	$$
Then,
\begin{enumerate}
\item If $T > T_0$,  system~\eqref{f11} is null controllable at time $T$.
\item If $T < T_0$,  system~\eqref{f11} is not null controllable at time $T$.
\end{enumerate}
\end{enumerate}
\end{theorem}

\begin{remark}
The proof of Theorem~\ref{t7} can be found in~\cite{O}. In fact, in that reference the author analyzes the controllability of $2 \times 2$ parabolic systems more general than system~\eqref{f11} with distributed or boundary scalar controls.
\end{remark}

\begin{remark}
The study of the controllability of system~\eqref{f11} is easier when $a \in L^2(0, \pi)$ does not satisfy condition~\eqref{f12}. In fact, we have the following property: system~\eqref{f11} is null controllable at time $ T >0$ if and only if the system is approximately controllable at this time, i.e., if and only if $ \Lambda_1 \cap \Lambda_2 = \emptyset$, with $\Lambda_1$ and $\Lambda_2$ given in the first item of Theorem~\ref{t7} (see~\cite{O}). In this case, we have that $T_0 = 0$ and the null controllability of the system for any $T > 0$. On the other hand, it is not difficult to check that we can apply Lemma~\ref{l2} to  the sequence $\Lambda$. As a consequence, the associated control cost $\mathcal{K} (T) $ for system~\eqref{f11} can be estimated as follows:
	$$
\exp \left( \frac{C_0} T  \right) \le \mathcal{K}(T) \le \exp \left( \frac{C_1} T  \right), \quad \forall T \in (0, \tau_0) ,
	$$
for appropriate positive constants $C_0$, $C_1$ and $\tau_0$ independent of $T$.
\end{remark}

As said before, our objective is to provide an estimate of the control cost of system~\eqref{f11} in a simple case ($a$ satisfying~\eqref{f12}). To this end, let us state a technical result whose proof can be found in~\cite{PT} (see also~\cite{O}):
\begin{lemma}\label{l3}
Given $\gamma \in (0, 1)$, there exists $a \in L^2(0, \pi)$ satisfying~\eqref{f12} such that 
	$$
\sigma(L) = \Lambda = \ens{k^2, k^2 + e^{- k^{2 \gamma} }}_{k \ge 1} .
	$$
\end{lemma}

From now on, we will take $\gamma \in (0,1) $ and will work with the function $a$ provided by the previous lemma. Observe that, in particular (see Theorem~\ref{t7}), $ \Lambda_1 \cap \Lambda_2 = \emptyset$ and the minimal time $T_0$ for system~\eqref{f11} associated to $a$ is
	$$
T_0 = \limsup \frac{- \log \left( e^{- k^{2 \gamma} } \right)}{k^2} = \lim \frac{k^{2 \gamma}}{k^2} = 0.
	$$
Thus, system~\eqref{f11} is approximately and null controllable at any time $T > 0$.

On the other hand, the real sequence $\Lambda$ ($\Lambda$ is given in Lemma~\ref{l3}) can be rearranged as an increasing sequence $\Lambda = \ens{\lambda_k}_{k \ge 1} \subset \R$ as follows:
	$$
\lambda_k = \left\{ 
	\begin{array}{ll}
\displaystyle \left( \frac{k + 1}2 \right)^2 & \hbox{if $k$ is odd}, \\
	\noalign{\smallskip}
\displaystyle \left( \frac{k }2 \right)^2 + e^{-(k/2)^{2 \gamma}} & \hbox{if $k$ is even}.
	\end{array}
	\right.
	$$

It is not difficult to check that $\Lambda$ satisfies conditions~\eqref{hyp 1}--\eqref{hyp 5},~\eqref{hyp 7} and~\eqref{hyp_6} in Theorem~\ref{t5}, with $q = 2$. Therefore, the previous result can be applied to $\Lambda $ obtaining the existence of a constant $C >0$ such that for any $T  > 0$ there exists a biorthogonal family $\ens{q_k}_{k \ge 1} $ to $\ens{e^{- \lambda_ k t}}_{k \ge 1}$ in $L^2(0,T)$ such that
	$$
\| q_k \|_{L^2(0,T)} \le C e^{C  \lceil k/2 \rceil } e^{C /T} e^{\lceil k/2 \rceil^{2 \gamma}} , \quad \forall k\geq 1,
	$$
where $\lceil \cdot \rceil$ is the ceiling function\footnote{$ \lceil x \rceil $ is the least integer greater than or equal to 
$x$.} (see~\cite{GBO} for the details).

Combining the previous estimate and the moment method, we deduce:
\begin{theorem}\label{t8}
Let us fix $\gamma \in (0,1) $ and take the function $a \in L^2(0, \pi)$ provided by Lemma~\ref{l3}. If we denote $\mathcal{K} (T) $ the control cost of system~\eqref{f11} in $L^2(0,T)$, then, there exist two positive constants $\tau_0$ and $C_1$ such that
	\begin{equation}\label{f13}
\mathcal{K}(T) \le \exp \left( \frac{C_1} T + \frac{C_1}{T^{\frac\gamma{1 - \gamma}}} \right), \quad \forall T \in (0, \tau_0) .
	\end{equation}
\end{theorem}

The proof of the previous result can be found in~\cite{GBO}. In fact, the previous estimate for $\mathcal{K} (T)$ is optimal in the following sense (see~\cite{GBO}):

\begin{theorem}\label{t9}
Under the assumptions of Theorem~\ref{t8}, there exists a positive constant $C_0$ such that
	\begin{equation}\label{f14}
\mathcal{K}(T) \ge \exp \left( \frac{C_0} T + \frac{C_0}{T^{\frac\gamma{1 - \gamma}}} \right), \quad \forall T \in (0, \tau_0) .
	\end{equation}
\end{theorem}

Again, a proof of this result can be found in~\cite{GBO}.

\begin{remark}
Observe that inequalities~\eqref{f13} and~\eqref{f14} are valid when $\gamma \in (0,1)$. In fact, we can write
\begin{enumerate}
\item If $\gamma \in (0, 1/2]$ one deduces the existence of positive constants $C_0, C_1$ such that
	$$
\exp \left( \frac{C_0} T  \right) \le \mathcal{K}(T) \le \exp \left( \frac{C_1} T  \right), \quad \forall T \in (0, \tau_0) ,
	$$
obtaining estimates for the control cost of system~\eqref{f11} similar to the heat equation (see Theorems~\ref{t2} and~\ref{t3}).
\item If $\gamma \in (1/2, 1)$, then, for new positive constants $C_0, C_1$ we have
	$$
\exp \left( \frac{C_0}{T^{\frac\gamma{1 - \gamma}}}  \right) \le \mathcal{K}(T) \le \exp \left( \frac{C_1}{T^{\frac\gamma{1 - \gamma}}}  \right), \quad \forall T \in (0, \tau_0) .
	$$
The previous expression proves that the control cost blows up when $\gamma \to 1^-$. This is natural because the minimal time for system~\eqref{f11} when $\gamma = 1$ is $T_0 = 1$ and the system is not null controllable at time $T$ when $T < 1$.
\end{enumerate}
\end{remark}

\begin{remark}
In this section we have analyzed the control cost for system~\eqref{f11} when $a \in L^2(0, \pi)$ is such that the sequence of  eigenvalues of the corresponding generator is given by $\Lambda = \ens{k^2, k^2 + e^{k^{2 \gamma} }}_{k \ge 1}$ with $\gamma \in (0,1)$. With this choice, the minimal time of null controllability of the system is $T_0 = 0$. In~\cite{GBO}, the authors analyze the general problem $a \in L^2(0, \pi)$ satisfying~\eqref{f12} which includes the case in which $T_0 >0$. The estimates for $\mathcal{K} (T)$ that can be obtained when $T_0 > 0$ are very close to those obtained in~\cite{L} for the heat equation (see~\cite{GBO}).
\end{remark}


\end{document}